\documentclass[11pt, twoside]{amsart}

\usepackage[utf8]{inputenc}
\usepackage[T1]{fontenc}
\usepackage{amssymb,amsmath,amstext, graphicx, hyperref}

\newtheorem{theor}{Theorem}
\newtheorem{lem}[theor]{Lemma}
\newtheorem{defin}[theor]{Definition}

\newtheorem{prop}[theor]{Proposition} 

\newtheorem{notation}[theor]{Notation}

\newtheorem{cor}[theor]{Corollary}
\newtheorem{rem}[theor]{Remark}

\newtheorem{fact}[theor]{Fact}

\newcommand{\dom}{\mathrm{dom}}
\newcommand{\acl}{\mathrm{acl}}
\newcommand{\dcl}{\mathrm{dcl}}
\newcommand{\tp}{\mathrm{tp}}

\newcommand{\es}{\emptyset}

\newcommand{\su}{\mathrm{SU}}

\newcommand{\nts}{\negthickspace}
\newcommand{\uhrc}{\nts \upharpoonright \nts}

\newcommand{\meq}{^{\mathrm{eq}}}

\newcommand{\mcM}{\mathcal{M}}
\newcommand{\mcN}{\mathcal{N}}

\newcommand{\mcP}{\mathcal{P}}

\newcommand{\us}{\underset}
\newcommand{\ind}{\raisebox{-2pt}[5pt][0pt]{$\smile$} \hspace*{-6.8pt}\raisebox{3pt}[5pt][0pt]{$|$} \; \: }
\newcommand{\nind}{\raisebox{-2pt}[5pt][0pt]{$\smile$} 
\hspace*{-6.8pt}\raisebox{3pt}[5pt][0pt]{$|$}\hspace*{-6.8pt}
\raisebox{3pt}[5pt][0pt]{$\diagup$} }

\newcommand{\rng}{\mathrm{rng}}

\title[Binary simple homogeneous structures]
{Binary simple homogeneous structures are supersimple with finite rank}

\author{Vera Koponen}

\address{Vera Koponen, Department of Mathematics, Uppsala University, Box 480,
75106 Uppsala, Sweden.}

\email{vera@math.uu.se}

\begin{document}

\begin{abstract}

Suppose that $\mcM$ is an infinite structure with finite relational vocabulary such that every relation symbol has arity at most 2.
If $\mcM$ is simple and homogeneous then its complete theory is supersimple with finite SU-rank which
cannot exceed the number of complete 2-types over the empty set. 
\end{abstract}

\maketitle

\section{Introduction}\label{Introduction}

\noindent
A first-order structure $\mcM$ will be called {\em homogeneous} (sometimes called {\em finitely homogeneous} or
{\em ultrahomogeneous}) if it is countable, has a finite vocabulary (signature) with only relation symbols (a {\em relational vocabulary})
and every isomorphism between finite substructures of $\mcM$ can be extended to an automorphism of $\mcM$. 
Although being countable and having a finite relational vocabulary is part of being homogeneous according to this definition,
we will sometimes repeat these assumptions.
If the vocabulary of $\mcM$ has only relations symbols that are unary or binary, then 
say that the vocabulary is binary and call $\mcM$ a {\em binary structure}.
For countable $\mcM$ with finite relational vocabulary, 
$\mcM$ is homogeneous if and only if $\mcM$ has elimination of quantifiers
\cite[Corollary 7.4.2]{Hod}. 
Via quantifier elimination one can see that infinite
homogeneous structures are $\omega$-categorical \cite{Hod}.
Moreover, a structure is homogeneous if and only if it is the
so called Fra\"{i}ss\'{e} limit of an ``amalgamation class'' of finite structures \cite{Fra, Hod}. 
Besides being interesting objects from a model theoretic point of view, homogeneous structures have been studied
in connection to Ramsey theory, constraint satisfaction problems,
permutation groups and topological dynamics. See \cite{BP, Che98, HN, Mac10, Nes} for surveys of homogeneous structures
and their applications.

We are far from a good understanding of homogeneous structures in general, although some particular classes of homogeneous
structures have been classified or are very well understood \cite{Che98, Gar, GK, JTS, Lach84, Lach97, LT, LW, Schm, Shee}.
The framework of model theoretic stability theory, later generalized to simplicity theory, 
gives tools which makes it possible to understand structures in a quite general context.
We say that an infinite structure is stable/simple if its complete first-order theory is stable/simple.
Lachlan and his collaborators used tools available for stable structures to work out a very detailed understanding
of infinite stable homogeneous structures; see for example the survey \cite{Lach97}.\footnote{
They considered every finite homogeneous structure to be stable 
so also finite homogeneous structures were part of their analysis, but for the finite ones the theory is
not as conclusive as for the infinite ones.} When saying that a structure is stable or simple we will from now on
assume that it is infinite.
The present work and \cite{AL, AK, Kop14} can be seen as a continuation and (to the extent possible)
generalization of the work on stable homogeneous structures.
This seems worthwhile since, on the one hand, stability/simplicity theoretic ideas appear to be useful
beyond the context of stable homogeneous structures, and, on the other hand, because
new phenomena arise in unstable simple homogeneous structures, and these new phenomena show that the
class of simple homogeneous structures is, in interesting ways, richer than the class of stable homogeneous structures.
This may be of interest to applications of homogeneous structures.
Some of these differences are discussed below.

All stable homogeneous structures are in fact $\omega$-stable and hence superstable.
This fact follows fairly quickly from the characterizations of these notions by counting types and the fact that homogeneous
structures have elimination of quantifiers. 
Somewhat more precisely, if $\mcM$ is (infinite) homogeneous and not $\omega$-stable, then, for
some countable set $A$ the set of  complete 1-types over $A$ is uncountable, and by elimination of quantifiers 
with respect to a finite relational vocabulary there
must be an atomic formula $\varphi$ such that there are uncountably many pairwise inconsistent
1-types over $A$ which use {\em only} the formula
$\varphi$. Shelah's ``unstable formula theorem'' \cite[Ch. II,Theorem 2.2]{She} now implies that 
$\mcM$ is not stable.
This argument cannot be generalized to prove that every (binary) simple homogeneous structure is supersimple,
because, by a well known characterization of stability, 
every {\em un}stable (first-order) theory has $2^\lambda$ complete
types over some set of parameters of cardinality $\lambda$, for every choice of infinite cardinal $\lambda$.

The following is essential for the theory of stable homogeneous structures, where $rk$ is
the rank used by Lachlan which is derived from Shelah's ``CP$( \ , 2)$-rank'': if $\mcM$ is stable and homogeneous then
(a) $rk(\mcM)$ is finite, and (b) every set with U-rank 1 and without a nontrivial definable equivalence relation is indiscernible.\footnote{The connection with Lachlan's terminology is the following: if the structure under consideration is stable and homogeneous, then
any set which has U-rank 1 and no nontrivial definable equivalence relation is {\em strictly minimal} (and hence indiscernible).}
Neither~(a) nor~(b) is true in general for (even binary) simple homogeneous structures.
For example, if $\mcM$ is the Rado graph, in model theory often called the {\em random graph}, 
then its universe is a set with U-rank 1 and without a nontrivial definable equivalence relation, but it is not indiscernible.
This failure can at least partially be blamed on the failure of simple unstable structures to have
``unique nondividing extensions of stationary types''; a property which all stable structures have
\cite[Ch. III, Corollary 2.9]{She}.
The failure of~(a) for unstable infinite structures is tightly connected to Shelah's ``unstable formula theorem'';
see Theorem 2.2, Theorem 3.2, Definition 3.4 and Exercise 3.8 in Chapter~II of~\cite{She}.

In spite of the failures of~(a) and~(b) for unstable simple homogeneous structures, these 
statements seem to point in the right direction. For example, Theorem~5.1 in \cite{AK} may be seen as a
version of~(b) in the case of binary simple homogeneous structures.
It also makes sense to use some notion of rank when studying simple homogeneous structures,
as will be further discussed below.
The so-called U-{\em rank} is important in many studies of stable infinite structures
and since for all $\omega$-categorical superstable $\mcM$ we have 
$\mathrm{U} (\mcM) \leq rk(\mcM)$ \footnote{
This inequality can be understood as follows.
By \cite[Corollary 6.48]{Pil}, U-rank coincides with Morley rank \cite[Definition 6.16]{Pil} if the structure
is $\omega$-categorical and superstable. Moreover, it is straightforward to see that $rk(\mcM)$ is at least
as big as the Morley rank of $\mcM$ (where the latter is the supremum of the Morley ranks of all 1-types consistent with
the complete theory of $\mcM$). Since every stable homogeneous structure $\mcM$ is $\omega$-categorical and
superstable, we get $\mathrm{U}(\mcM) \leq rk(\mcM)$.
}
it follows from~(a) that every infinite stable homogeneous
structure $\mcM$ has finite U-rank.
The U-rank also makes sense for simple structures, but in the context of simple structures it is
usually called SU-{\em rank}, a convention which we follow here.
If $T$ is a simple theory then the SU-rank of a complete type (over any set of parameters) is an
ordinal or $\infty$, where $\infty$ is understood to be larger than every ordinal \cite{Cas, Wag}.
If every type has ordinal valued SU-rank then the theory is called {\em supersimple}. 
The SU-rank of a supersimple theory $T$ is the supremum of the SU-ranks of all 1-types over $\es$ 
(that are realized by ``real elements'').
Experience has shown that properties of a simple theory with finite SU-rank can often be analysed
via properties of types of SU-rank 1.
For example, an $\omega$-categorical simple theory $T$ with finite SU-rank is 1-based if and only if
all types of SU-rank 1 are 1-based (sometimes called modular) \cite[Corollary 4.7]{HKP}.

All known simple homogeneous structures have complete theories which are supersimple and have finite SU-rank (and are 1-based). 
The most famous example is probably the random graph, which has SU-rank 1.
Aranda L\'{o}pes \cite[Theorem 3.2.7]{AL} has proved that if $\mcM$ is binary, homogeneous and
supersimple, then the SU-rank of the complete theory of $\mcM$ cannot be $\omega^\alpha$ for any ordinal $\alpha \geq 1$.
I am not aware of other results in this direction for simple (unstable) homogeneous structures.
The main result of this article is the following, where $S_2(T)$ is the set
of complete 2-types over $\es$ with respect to the theory $T$:

\begin{theor}\label{binary homogeneous simple structures have finite rank}
Suppose that $\mcM$ is a countable, binary, homogeneous and simple structure.
Let $T$ be the complete theory of $\mcM$.
Then $T$ is supersimple with finite SU-rank which is at most $|S_2(T)|$.
\end{theor}

\noindent
Since homogeneous structures have elimination of quantifiers it follows from Theorem~1 that,
for every binary finite relational vocabulary $V$, every simple homogeneous $V$-structure is supersimple
with SU-rank at most $2^{c|V|}$, where $c$ is a constant that depends only on $V$.

With the above theorem at hand some questions about a binary homogeneous and simple structure $\mcM$
can be studied by asking the analogous questions for types of SU-rank 1.
Let $T$ be the complete theory of $\mcM$.
By a result of Hart, Kim and Pillay \cite{HKP}, $T$ is 1-based (called `modular' in \cite{HKP}) if and only if 
every type of SU-rank 1 is 1-based. Moreover, by also involving work of Macpherson \cite{Mac91} and De Piro and Kim \cite{PK}
it follows that $T$ is 1-based if and only if $T$ has trivial dependence if and only if 
every type of SU-rank 1 has trivial pregeometry (see for example \cite[Section 2.3]{Kop14} for definitions and more explanation).
If for every homogeneous simple structure $\mcM$ its complete theory has trivial dependence and finite SU-rank,
then the behavior of dependence in simple homogeneous structures parallels that of stable homogeneous structures
(see \cite{Lach97} for a survey of stable homogeneous structures).
The reader is referred to \cite{AK, AL, Kop14} for more results about simple homogeneous structures.

It is natural to ask whether the `binarity' assumption in 
Theorem~\ref{binary homogeneous simple structures have finite rank}
is necessary,
especially as Cherlin and Lachlan proved that $rk(\mcM)$, and hence $\mathrm{U}(\mcM)$, is finite
for every stable homogeneous $\mcM$ \cite{CL}.
Their proof relies heavily on the classification of finite simple groups and on the possibility, in their context,
to reduce certain problems to questions about the automorphism group of a finite structure.
(However, in the binary case this classification is not needed \cite{LS}.)
Moreover, the rank $rk$ that they consider makes sense for finite and infinite 
structures while the definition of SU-rank (and U-rank) presupposes
that the structure in question is infinite.
In addition, the reduction to finite structures in \cite{CL} may 
depend on a property called  {\em smooth approximability (by finite substructures)}
\cite{CH, KLM} which holds for all stable homogeneous structures, but not for all simple homogeneous structures
(the random bipartite graph is not smoothly approximable \cite[p 457]{KLM}).
So there are some seemingly difficult obstacles if one tries to generalize 
Theorem~\ref{binary homogeneous simple structures have finite rank}
to nonbinary structures via the approach of \cite{CL}.
On the other hand, the more advanced state of stability/simplicity theory at present compared to the 1980'ies
may offer tools with which one can bypass or mitigate these obstacles.
The proofs of this article and of \cite{AL, AK, Kop14} show that simplicity
(of structures) and binarity have strong consequences when combined (via application of the  ``independence theorem''
of simple structures \cite[Theorem 2.5.20]{Wag}); these consequences
do not in an obvious way transfer to the context of nonbinary structures.\footnote{
At least one result in the binary case, namely Theorem~5.1 in \cite{AK}, cannot be ``naturally translated''
to the nonbinary case,
as witnessed by the ``generic pyramid-free 3-hypergraph'' which is simple with SU-rank 1;
see \cite[Section 3]{Djo} for details.}
But nevertheless I tend to believe that the assumption of binary vocabulary in
Theorem~\ref{binary homogeneous simple structures have finite rank}
is not necessary.

The structure of this article is as follows. 
The next section recalls the necessary background and explains some notation and terminology.
Theorem~\ref{binary homogeneous simple structures have finite rank}
is proved in Section~\ref{Finiteness of rank}. 
In Section~\ref{An auxilliary result about independence} a technical result about independence in
simple homogeneous structures is proved, which is then used in the proof of 
Theorem~\ref{binary homogeneous simple structures have finite rank}.

\section{Preliminaries}\label{Preliminaries}

\noindent
The notation and terminology that we use is more or less standard, but nevertheless we explain some notational issues here.
Structures are denoted by $\mcM$ or $\mcN$ (or $\mcM\meq$ or $\mcN\meq$ if we deal with imaginaries) 
and their universes are denoted by $M$ or $N$ (or $M\meq$ or $N\meq$). Sometimes we attach indices to the letters.
Finite sequences (tuples) of elements are denoted $\bar{a}, \bar{b}, \ldots$ (and finite sequences of variables $\bar{x}, \bar{y}, \ldots$)
while $a, b, \ldots$ denote elements from some structure.
By `$\bar{a} \in A$' we mean that all elements in the sequence $\bar{a}$ belong to the set $A$.
If we want to show that the length of $\bar{a}$ is $n$ then we may write $\bar{a} \in A^n$.
For a set $A$, $|A|$ is its cardinality and for a sequence $\bar{a}$, $|\bar{a}|$ is its length.
For a sequence $\bar{a}$, $\rng(\bar{a})$ denotes the set of elements occuring in the sequence.
The {\em maximal arity} of a finite relational vocabulary is, of course, the maximum of the arities of
the relation symbols in the vocabulary.

Suppose that $\mcM$ is a structure, $\bar{a} \in M$ and $A \subseteq M$.
Then $\acl_\mcM(A)$, $\dcl_\mcM(A)$ and $\tp_\mcM(\bar{a} / A)$ denote the 
{\em algebraic closure} of $A$ with respect to $\mcM$, the {\em definable closure} of $A$ with respect to $\mcM$ and
the {\em complete type of $\bar{a}$ over $A$} with respect to $\mcM$, respectively. 
By $tp_\mcM^{at}(\bar{a} / A)$ we mean the restriction of 
$\tp_\mcM(\bar{a} / A)$ to atomic formulas.
We often write $\tp_\mcM(\bar{a})$ instead of $\tp_\mcM(\bar{a} / \es)$ (and similarly for `$tp_\mcM^{at}$').
With $\mcM \uhrc A$ we denote the substructure of $\mcM$ which is generated by $A$.
Observe that if the vocabulary of $\mcM$ is relational,
then $\tp_\mcM^{at}(a_1, \ldots, a_n) = \tp_\mcM^{at}(b_1, \ldots, b_n)$ is equivalent to 
saying that the map $a_i \mapsto b_i$ is an isomorphism from $\mcM \uhrc \{a_1, \ldots, a_n\}$
to $\mcM \uhrc \{b_1, \ldots, b_n\}$.

Still assume that $A \subseteq M$.
By $S_n^\mcM(A)$ we denote the set of all complete $n$-types over $A$ which are realized in
some elementary extension of $\mcM$.
For a complete theory $T$ we let $S_n(T)$ be the set of all complete $n$-types (without parameters) of $T$.
This means that if $\mcM \models T$, then $S_n(T) = S_n^{\mcM}(\es)$.
Suppose that $R \subseteq M^k$. Then we say that
$R$ is {\em $A$-definable} (in $\mcM$) if there are a formula $\varphi(\bar{x}, \bar{y})$ without parameters
and $\bar{a} \in A$ such that $R = \{\bar{b} \in M^k : \mcM \models \varphi(\bar{b}, \bar{a})\}$.
A structure $\mcM$ is called {\em $\omega$-categorical}, {\em simple} or {\em supersimple}, respectively, if its complete
theory, denoted $Th(\mcM)$, has that property.

We refer to \cite{Hod} (for example) for unexplained basic notions and notation of model theory,
and to \cite{Cas, Wag} for basic concepts and results from simplicity theory.

We have to work a little bit with imaginary elements, in order to show that we can avoid them in the crucial 
part of the proof of Theorem~\ref{binary homogeneous simple structures have finite rank}.
As usual $\mcM\meq$ denotes the extension of $\mcM$ by imaginary elements.
The approach to imaginary elements that we adopt is that of \cite{Hod, She} in which we do {\em not}
introduce variables of different sorts but instead use unary predicates to ``point out'' the different sorts.
This approach is also used in \cite{AK} where it is explained in more detail. 
The following fact which we will use is also explained in some more detail in \cite{AK}.

\begin{fact}\label{facts about simple omega-categorical structures}
Suppose that $\mcM$ is $\omega$-categorical. Then:\\
(i) For all $\bar{a}, \bar{b} \in M$, 
$\tp_\mcM(\bar{a}) = \tp_\mcM(\bar{b})$ if and only if $\tp_{\mcM\meq}(\bar{a}) = \tp_{\mcM\meq}(\bar{b})$
($\omega$-categoricity is not needed for this part). \\
(ii) If $B \subseteq M\meq$ is finite and $\bar{a} \in M\meq$,
then $\tp_{\mcM\meq}(\bar{a} / \acl_{\mcM\meq}(B))$ is isolated. \\
(iii) If $B \subseteq \mcM\meq$ is finite, $n < \omega$ and $p \in S_n^{\mcM\meq}(\acl_{\mcM\meq}(B))$ is realized 
in $\mcN\meq$ for some $\mcN \succcurlyeq \mcM$,
then $p$ is realized in $\mcM\meq$.
\end{fact}

\noindent
Let $T$ be a simple theory. 
For every $\mcM \models T$, $A \subseteq M\meq$ and $p \in S_n^{\mcM\meq}(A)$, 
there is a notion of {\em $\su$-rank} of $p$, denoted $\su(p)$
(a definition is found in \cite{Cas, Wag}). 
We abbreviate $\su(\tp_{\mcM\meq}(\bar{a} / A))$
with $\su(\bar{a} / A)$.
For every type $p$,  $\su(p)$ is either ordinal valued or undefined (or alternatively given the value $\infty$).
$T$ is {\em supersimple} if and only if for every $0 < n < \omega$ and 
$p \in S_n(T)$, $\su(p)$ is ordinal valued (by \cite[Proposition~13.13]{Cas} or \cite[Theorem~5.1.5]{Wag}, 
and the facts that $p \subseteq q$ implies that $\su(p) \geq \su(q)$ and if $\mcM \models T$ and
$\bar{a} \in M\meq$, then there is
$\bar{a}' \in M$ such that $\bar{a} \in \dcl_{\mcM\meq}(\bar{a}')$ and hence 
$\su(\bar{a}) \leq \su(\bar{a}')$).

The {\em $\su$-rank of $T$} is the supremum of $\{\su(p) : p \in S_1(T)\}$. 
If the SU-rank of $T$ is finite then it follows from the Lascar inequalities \cite{Cas, Wag} that $\su(p)$ is finite
for every $p \in S_n(T)$ and every $n < \omega$;
so in particular, $T$ is supersimple.

\section{An auxilliary result about independence}\label{An auxilliary result about independence}

\noindent
In this section we prove a result (generalizing \cite[Theorem 3.3]{AK} and its proof) which will be used 
in the proof of the main theorem. 
Actually we will only use its corollary to binary structures, but nevertheless prove the more general version since it
may be useful in the future.
A slightly weaker version of Corollary~\ref{pairwise independence implies independence}
has been proved earlier by Aranda Lopez \cite{AL}.

We consider a generalization of the independence theorem for simple theories, namely the 
`strong $n$-dimensional amalgamation property for Lascar strong types', 
studied by Kolesnikov in \cite[Definition 4.3]{Kol05}. 
In the present context of homogeneous structures,
as distinct from that of Kolesnikov,  a `Lascar strong type' corresponds to a `type over an algebraically closed set'.
The notation $\mcP(S)$ denotes the powerset of $S$, and we let $\mcP^-(S) = \mcP(S) - \{S\}$. 
Every $n < \omega$ is identified with the set $\{0, \ldots, n-1\}$, so the notation $\mcP(n)$ makes sense. 
For a type $p$, $\dom(p)$ denotes the set of all parameters that occur in formulas in $p$.

\begin{defin}{\rm
Let $T$ be an $\omega$-categorical  and simple complete theory and let $n < \omega$. \\
(i) A set of types $\{p_s(\bar x) | s\in \mcP^-(n)\}$ (with respect to $\mcM\meq$ for some $\mcM \models T$) 
is called an {\em $n$-independent system of strong types over $A$} (where $A \subseteq \mcM\meq$)
if it satisfies the following properties:
\begin{itemize}
\item[(a)] $\dom(p_\emptyset) = A$ and $\dom(p_s)$ is algebraically closed in $\mcM\meq$ for every $s \in \mcP^-(n)$.
\item[(b)] for all $s,t \in \mcP^-(n)$ such that $s\subseteq t$, $p_t$ is a nondividing extension of $p_s$.
\item[(c)] for all $s,t \in \mcP^-(n)$, $\dom(p_s)\us{\dom(p_{s\cap t})}{\ind} \dom(p_t)$.
\item[(d)] for all $s,t \in \mcP^-(n)$, $p_s$ and $p_t$ extend the same type over $\acl_{\mcM\meq}(\dom(p_{s\cap t}))$.
\end{itemize}
(ii) We say that $T$ (and any $\mcN \models T$) has the {\em  n-dimensional amalgamation property for strong types} 
if for every $\mcM \models T$ and every
$n$-independent system of strong types $\{p_s(\bar x) | s\in \mcP^-(n)\}$ over some set $A \subseteq M\meq$,
there is a type $p^*$ which is a nondividing extension of $p_s$ for each $s\in P^-(n)$.
}\end{defin}

\begin{rem}\label{remark that simplicity implies the 2-dimensional amalgamation property}{\rm
The independence theorem (in its general setting when the sets of parameters of the given types may be infinite \cite{Cas, Wag})
implies that every $\omega$-categorical and simple theory has the 2-dimensional amalgamation property for strong types.
(This relies on the fact that since $\omega$-categorical theories have elimination of hyperimaginaries \cite{Cas, Wag} we can
replace the `bounded closure' with `algebraic closure'.)
}\end{rem}

\begin{prop}\label{pregeometries are degenerate if the structure is homogeneous}
Suppose that $\mcM$ has a finite relational vocabulary with maximal arity~$\rho$.
Also assume that $\mcM$ is countable, homogeneous and simple and has the 
$\rho$-dimensional amalgamation property for strong types.
Let $0 < n < \omega$, $\bar{a}_0, \ldots, \bar{a}_n \in M$ and suppose that for every 
$s \subseteq \{0, \ldots, n\}$ such that $|s| \leq \rho$,
$\{\bar{a}_i : i \in s\}$ is independent over $B \subseteq M$ and that $\rng(\bar{a}_i) \cap \rng(\bar{a}_j) = \es$
whenever $i < j \leq n$.
Then $\{\bar{a}_0, \ldots, \bar{a}_n\}$ is independent over $B$.
\end{prop}

\noindent
{\bf Proof.}
Suppose that $\mcM$, $\rho$ and $\bar{a}_0, \ldots, \bar{a}_n \in M$ satisfy the assumptions of the proposition.
Recall that $\rho \geq 2$ by our definition in Section~\ref{Preliminaries}.
We use induction on $n$.
The base case is when $n < \rho$ and then the conclusion is evident.
So suppose that $n \geq \rho$.
By the induction hypothesis, every proper subset of $\{\bar{a}_0, \ldots, \bar{a}_n\}$ is independent over $B$.
For a contradiction suppose that $\{\bar{a}_0, \ldots, \bar{a}_n\}$ is not independent over $B$.
Then for some $i \leq n$, $\bar{a}_i \underset{B}{\nind} \{\bar{a}_j : j \leq n \text{ and } j \neq i\}$.
Without loss of generality assume that $i = n$, so
\begin{equation}\label{a-n is not independent from the other b-i}
\bar{a}_n \underset{B}{\nind} \{\bar{a}_i : i < n\}
\end{equation}
and, by the induction hypothesis, 
\begin{equation}\label{all b-i for i < n is an independent set over B}
\{ \bar{a}_i : i < n \} \text{ is independent over } B. 
\end{equation}
The induction hypothesis also implies that 
\begin{equation}\label{a-n does not fork any at most n-1 of the other tuples}
\text{for all } s \in \mcP^-(n), \ \bar{a}_n \underset{B}{\ind} \{ \bar{a}_i : i \in s\}. 
\end{equation}

For each $s \in \mcP^-(\rho)$, let 
\[A_s = \acl_{\mcM\meq}\big(\{\bar{a}_i : i \in s\} \cup  \{\bar{a}_\rho, \ldots, \bar{a}_{n-1}\} \cup B\big).\]

\noindent
{\bf Claim 1.} For all $s, t \in \mcP^-(\rho)$, $A_s \underset{A_{s \cap t}}{\ind} A_t$.

\medskip

\noindent
{\em Proof of Claim~1.}
Suppose that $s, t \in \mcP^-(\rho)$ and $A_s \underset{A_{s \cap t}}{\nind} A_t$.
Then $t \setminus s \neq \es$.
If $|t \setminus s| > 1$ then there is $t' \subset t$ such that
$t' \cap s = t \cap s$, $|t \setminus t'| = 1$ and $|t' \setminus s| > 0$.
By transitivity of dividing,
$A_s \underset{A_{s \cap t'}}{\nind} A_{t'}$ or 
$A_s \underset{A_{t'}}{\nind} A_t$;
in the latter case $A_{s \cup t'} \underset{A_{(s \cup t') \cap t}}{\nind} A_t$
(because $(s \cup t') \cap t = t'$).
In the first case, $|t' \setminus s| < |t \setminus s|$.
In the second case, $|t \setminus (s \cup t')| < |t \setminus s|$.
By induction on $|t \setminus s|$ we therefore find $s', t' \in \mcP^-(\rho)$ such that
$A_{s'} \underset{A_{s' \cap t'}}{\nind} A_{t'}$ and $|t' \setminus s'| = 1$.
By monotonicity of dividing, $A_{s'} \underset{B}{\nind} A_{t'}$ where 
$A_{t'} \setminus A_{s'}$ contains exactly one tuple from
$\{\bar{a}_i : i < \rho\}$
(because $|t' \setminus s'| = 1$).
Hence $\big\{\bar{a}_i : i \in s' \cup t' \cup \{\rho, \ldots, n-1\}\big\}$ 
is {\em not} independent over $B$, which contradicts~(\ref{all b-i for i < n is an independent set over B})
since $s' \cup t' \subseteq \rho \leq n$.
$\text{ }$ \hfill $\square$

\medskip

\noindent
Now we verify that 
\[ \big\{ \tp(\bar{a}_n / A_s) : s \in \mcP(\rho)^- \big\} \]
is a $\rho$-independent system of strong types over 
$\acl_{\mcM\meq}(B \cup \{\bar{a}_\rho, \ldots, \bar{a}_{n-1}\})$.
Properties~(a) and~(d) follow directly from the definition of $A_s$ for $s \in \mcP^-(\rho)$;
(b) follows from~(3) and monotonicity of dividing; and~(c) follows from Claim~1.
Since $\mcM$ has the $\rho$-dimensional amalgamation property for strong types 
there is $p^*(\bar{x})$ which is a nondividing extension of $\tp(\bar{a}_n / A_s)$ for each $s \in \mcP^-(\rho)$.
Without loss of generality we may assume that $\dom(p^*) = \bigcup_{s \in \mcP^-(\rho)} A_s$.
Then, by Fact~\ref{facts about simple omega-categorical structures}, we find
$\bar{a} \in M$ which realizes $p^*$. 
Since $p^*$ does not divide over $A_\es$ it follows that 
\[
\bar{a} \underset{B \cup \{\bar{a}_\rho, \ldots, \bar{a}_{n-1}\}}{\ind} \{\bar{a}_i : i < \rho\}.
\]
This together with~(\ref{a-n does not fork any at most n-1 of the other tuples})
and transitivity gives $\bar{a} \underset{B}{\ind} \{\bar{a}_i : i < n\}$.
To sum up, we have:
\begin{align}\label{a satisfying amalgamation}
&\tp(\bar{a} / A_s) \ = \ \tp(\bar{a}_n / A_s) \ \text{ for every $s \in \mcP^-(\rho)$, and} \\
&\bar{a} \underset{B}{\ind} \{\bar{a}_i : i < n \}. \nonumber
\end{align}
Let 
\begin{align*}
C &= B \cup \rng(\bar{a}_0) \cup \ldots \cup \rng(\bar{a}_{n-1}) \cup \rng(\bar{a}_n) \ \text{ and}\\
C' &= B \cup \rng(\bar{a}_0) \cup \ldots \cup \rng(\bar{a}_{n-1}) \cup \rng(\bar{a}).
\end{align*}
Recall that our assumptions imply that $\rng(\bar{a}_n) \cap \rng(\bar{a}_i) = \es$ for 
all $i < n$. From~(\ref{a satisfying amalgamation}) it follows that 
$\rng(\bar{a}) \cap \rng(\bar{a}_i) = \es$ for all $i < n$

\medskip

\noindent
{\bf Claim 2.} The bijection $f : \mcM \uhrc C \to \mcM \uhrc C'$ defined by $f(x) = x$ for all
$x \in B \cup \rng(\bar{a}_0) \cup \ldots \cup \rng(\bar{a}_{n-1})$ and $f(\bar{a}_n) = \bar{a}$
is an isomorphism.

\medskip

\noindent
{\em Proof of Claim~2.}
By assumption every relation symbol has arity at most $\rho$. 
If $\bar{b} \in M^k$ where $k \leq \rho$, then $\rng(\bar{b})$ can have nonempty intersection with at most
$\rho$ of the sets $\rng(\bar{a}_i)$ for $i \leq n$, and similarly if we replace $\bar{a}_n$ with $\bar{a}$.
Therefore~(\ref{a satisfying amalgamation}) implies that for every relation symbol $R$ of arity $k$ and $\bar{b} \in C^k$, 
$\mcM \uhrc C \models R(\bar{b})$ if and only if $\mcM \uhrc C' \models R(f(\bar{b}))$, so $f$ is an isomorphism.
\hfill $\square$

\bigskip

\noindent
Since $\mcM$ is homogeneous and $B$ is finite, there is an automorphism $g$ of $\mcM$ which extends $f$ from Claim~2.
Then  $g(\bar{a}_n) = \bar{a}$ and $g$ fixes 
$B \cup \rng(\bar{a}_0) \cup \ldots \cup \rng(\bar{a}_{n-1})$ pointwise.
However, since dividing is invariant under automorphisms, this contradicts~(\ref{a-n is not independent from the other b-i})
and the second part of~(\ref{a satisfying amalgamation}).
\hfill $\square$
\\

\noindent
By Proposition~\ref{pregeometries are degenerate if the structure is homogeneous}
and Remark~\ref{remark that simplicity implies the 2-dimensional amalgamation property}
we get:

\begin{cor}\label{pairwise independence implies independence} {\rm (a slight strengthening of \cite[Proposition 3.1.4]{AL})}
Suppose that $\mcM$ is binary, countable, homogeneous and simple.
Let $0 < n < \omega$, $\bar{a}_0, \ldots, \bar{a}_n \in M$, let $B \subseteq M$ be finite and suppose that for all $i < j \leq n$,
$\bar{a}_i \underset{B}{\ind} \bar{a}_j$ and $\rng(\bar{a}_i) \cap \rng(\bar{a}_j) = \es$.
Then $\{\bar{a}_0, \ldots, \bar{a}_n\}$ is independent over~$B$.
\end{cor}

\noindent
The (anyway natural) assumption in Proposition~\ref{pregeometries are degenerate if the structure is homogeneous}
and Corollary~\ref{pairwise independence implies independence} that 
if $i \neq j$ then $\rng(\bar{a}_i) \cap \rng(\bar{a}_j) = \es$ could be removed, but at the cost of 
complicating the argument a little bit.

\section{Finiteness of rank}\label{Finiteness of rank}

\noindent
In this section we prove the main result:

\medskip

\noindent
{\bf Theorem~\ref{binary homogeneous simple structures have finite rank}.}{\em \
Suppose that $\mcM$ is a countable, binary, homogeneous and simple structure.
Let $T = Th(\mcM)$.
Then $T$ is supersimple with finite SU-rank which is at most $|S_2(T)|$.
}

\medskip

\noindent
Note that the bound on the rank need not be sharp: For every $k < \omega$ there is a binary random structure 
$\mcM$ such that $|S_2(Th(\mcM))| > k$, but the SU-rank of $Th(\mcM)$ is~1.
(See \cite[Section 2.3]{AK} for the exact meaning of binary random structure.)

We give the proof of Theorem~\ref{binary homogeneous simple structures have finite rank} 
in Section~\ref{proof of the main theorem}. Before that we do some preparatory work, including
introducing the notion of `preweight' which has a crucial role in the proof of 
Theorem~\ref{binary homogeneous simple structures have finite rank},
more precisely in Lemma~\ref{nondividing relation becomes dividing relation}.

\begin{lem}\label{homogeneity and simplicity is preserved under adding definable relation symbols}
Suppose that $\mcM$ is a countable homogeneous and simple  $V$-structure.
Let $0 < n < \omega$ and suppose that $A \subseteq M^n$ is a $\es$-definable relation.
Let $\mcM'$ be the expansion of $\mcM$ to the vocabulary $V \cup \{R_A\}$ where $(R_A)^{\mcM'} = A$.
Then $\mcM'$ is homogeneous and simple.
\end{lem}

\noindent
{\bf Proof.}
It is well known that simplicity is preserved if one adds relation symbols which are interpreted as relations that are $\es$-definable
in the original language.
One way of seeing this is to consider the {\em tree property} which is equivalent to {\em not} being simple \cite{Cas, Wag}: 
If a $(V \cup \{R_A\}$)-formula $\varphi'$ has the
tree property with respect to $Th(\mcM')$, then the $V$-formula $\varphi$ obtained by replacing
every occurence of $R_A$ with the $V$-formula which defines $A$ has the tree property with respect to $Th(\mcM)$.

Now suppose that $\bar{a}, \bar{b} \in M$ and $\tp^{at}_{\mcM'}(\bar{a}) = \tp^{at}_{\mcM'}(\bar{b})$.
Then $\tp^{at}_\mcM(\bar{a}) = \tp^{at}_\mcM(\bar{b})$ 
and as $\mcM$ is homogeneous there is an automorphism $f$  
of $\mcM$ such that $f(\bar{a}) = \bar{b}$.
Since $A \subseteq M^n$ is $\es$-definable, $f$ preserves $A$ setwise. Since $(R_A)^{\mcM'} = A$ it follows that
$f$ is an automorphism of $\mcM'$, so $\mcM'$ is homogeneous.
\hfill $\square$

\begin{defin}\label{definition of preweight}{\rm
Suppose that $\mcM$ is a simple structure.
Let $\bar{a} \in M\meq$, $B \subseteq M\meq$ and suppose that $\kappa$ is a cardinal.
The {\em preweight of $\bar{a}$ over $B$ (with respect to $\mcM$)}, denoted $pw(\bar{a} / B)$, is at least $\kappa$ if
there are $\mcN \succcurlyeq \mcM$ and a sequence $(\bar{a}_i : i < \kappa)$ in $\mcN\meq$ which is independent over $B$ and such that
$\bar{a} \underset{B}{\nind} \bar{a}_i$ for all $i < \kappa$.
We write $pw(\bar{a} / B) = \kappa$ if $pw(\bar{a} / B) \geq \kappa$ and $pw(\bar{a} / B) \not\geq \kappa^+$
(where $\kappa^+$ is the least cardinal greater than $\kappa$).
}\end{defin}

\begin{lem}\label{the preweight is finite} {\rm (D. Palac\'{i}n \cite[Lemma 2.14]{Pal})}
Let $V$ be a countable vocabulary. 
Suppose that $\mcM$ is a $V$-structure which is $\omega$-categorical and simple.
Let $\bar{a} \in M\meq$ and let $B \subseteq M\meq$ be finite.
Then $pw(\bar{a} / B) < \omega$.
\end{lem}

\subsection{Proof of Theorem~\ref{binary homogeneous simple structures have finite rank}}\label{proof of the main theorem}

Let $\mcM$ be a countable, binary, simple and homogeneous structure and let $T = Th(\mcM)$.
Moreover, by Lemma~\ref{homogeneity and simplicity is preserved under adding definable relation symbols} 
and the fact that the SU-rank of $T$ only depends on
which relations in models of $T$ are $\es$-definable (because of the definition of dividing), 
it follows that we may assume that 
\begin{align}\label{every 2-type is isolated by an atomic formula}
&\text{for every $p(x,y) \in S_2(T)$ there is a binary relation symbol $R_p$ such that} \\
&\text{$p(x,y)$ is isolated by $R_p(x,y)$.} \nonumber
\end{align}

\noindent
Let 
\[ t \ = \ |S_2(T)|. \]

\noindent
To prove that the SU-rank of $T$ is at most $t$ we need to prove (by the finite character of dividing/forking)
that there do not exist 
$\mcN \models T$, $a \in N$ and {\em finite} sets $\es = B_0 \subset B_1 \subset \ldots \subset B_{t+1} \subseteq N\meq$ such 
that $\tp_{\mcN\meq}(a / B_{n+1})$ divides over $B_n$ for every $n < t+1$.
As explained in the end of Section~\ref{Preliminaries}, supersimplicity follows from this.
The first step in the proof is to show that it suffices to consider the case
when $\mcN = \mcM$ and $B_n \subseteq M$ for all $n \leq t+1$.
This is taken care of by Lemmas~\ref{no long dividing sequence}
and~\ref{no long dividing sequence in M}.

\begin{lem}\label{no long dividing sequence}
Suppose that there are $\mcN \models T$, 
$a \in N$ and finite subsets 
\[ \es = B_0 \subset B_1 \subset \ldots \subset B_{t+1} \subset N\meq \] 
such that 
$\tp_{\mcN\meq}(a / B_{n+1})$ divides over $B_n$ for every $n < t+1$.
Then there are $a' \in M$ and finite 
\[ \es = B'_0 \subset B'_1 \subset \ldots \subset B'_{t+1} \subset M\meq \] 
such that 
$\tp_{\mcM\meq}(a' / B'_{n+1})$ divides over $B'_n$ for every $n < t+1$.
\end{lem}

\noindent
{\bf Proof.}
Suppose that $\mcN \models T$, $a \in N$ and that $B_n$, $n \leq t+1$ satisfy the assumptions of the lemma.
Without loss of generality we may assume that $\mcM \preccurlyeq \mcN$ (and $\mcM\meq \preccurlyeq \mcN\meq$).
Let $\bar{b}_{t+1}$ enumerate $B_{t+1}$. 
Then $\tp_{\mcN\meq}(a, \bar{b}_{t+1}) \in S_k^{\mcN\meq}(\es) = S_k^{\mcM\meq}(\es)$
for suitably chosen $k$, so by Fact~\ref{facts about simple omega-categorical structures},
$\tp_{\mcN\meq}(a, \bar{b}_{t+1})$ is also realized in $\mcM\meq$ by some $a', \bar{b}'_{t+1} \in M\meq$.
Then $\rng(\bar{b}'_{t+1})$ contains sets $B'_n$, for $n \leq t+1$, which satisfy the conclusion of the lemma.
$\text{  }$ \hfill $\square$

\begin{lem}\label{no long dividing sequence in M}
Suppose that there are $a \in M$ and finite subsets 
\[ \es = B_0 \subset B_1 \subset \ldots \subset B_{t+1} \subset M\meq \] 
such that 
$\tp_{\mcM\meq}(a / B_{n+1})$ divides over $B_n$ for every $n < t+1$.
Then there are finite subsets 
\[ \es = B'_0 \subset B'_1 \subset \ldots \subset B'_{t+1} \subset M \] 
such that
$\tp_\mcM(a / B'_{n+1})$ divides over $B'_n$ for every $n < t+1$.
\end{lem}

\noindent
{\bf Proof.}
Suppose that there are $a \in M$ and finite subsets 
$\es = B_0 \subset B_1 \subset \ldots \subset B_{t+1} \subset M\meq$ such that
$\tp_{\mcM\meq}(a / B_{n+1})$ divides over $B_n$ for every $n < t+1$.
For every $n \leq t+1$ there is finite $B'_n \subset M$ such that $B_n \subseteq \dcl_{\mcM\meq}(B'_n)$.
By enumerating $B'_n$ as $\bar{b}'_n$ and using the existence of a nondividing 
extension of $\tp_{\mcM\meq}(\bar{b}'_n/ B_n)$ to $B_n \cup \{a\}$ 
(and Fact~\ref{facts about simple omega-categorical structures}) we may also assume
that $a \underset{B_n}{\ind} B'_n$ for all $n \leq t+1$.
Suppose for a contradiction that $a \underset{B'_n}{\ind} B'_{n+1}$ for some $n < t+1$.
Then $a \underset{B_nB'_n}{\ind} B'_{n+1}$ because $B_n \subseteq \dcl_{\mcM\meq}(B'_n)$.
By transitivity of dividing 
(and since $a \underset{B_n}{\ind} B'_n$)
we get $a \underset{B_n}{\ind} B'_nB'_{n+1}$, so
$a \underset{B_n}{\ind} B'_{n+1}$.
As $B_{n+1} \subseteq \dcl_{\mcM\meq}(B'_{n+1})$
we have $a \underset{B_n}{\ind} B_{n+1}$ which contradicts the assumption.
Hence $a \underset{B'_n}{\nind} B'_{n+1}$ for every $n < t+1$.
Note that the construction does not guarantee that $B'_n \subseteq B'_{n+1}$. But by letting
$B''_n = \bigcup_{i \leq n}B'_n$ we get $B''_n \subseteq B''_{n+1}$ and still have $a \underset{B''_n}{\nind} B''_{n+1}$
for all $n < t+1$.
\hfill $\square$
\\

\noindent
From Lemmas~\ref{no long dividing sequence}
and~\ref{no long dividing sequence in M}
it follows that to prove that the SU-rank of $T = Th(\mcM)$ is at most $t$ 
it suffices to prove:

\begin{lem}\label{no long dividing sequence of real elements in M}
There do not exist $a \in M$ and finite subsets 
\[ \es = B_0 \subset B_1 \subset \ldots \subset B_{t+1} \subset M \] 
such that 
$\tp_\mcM(a / B_{n+1})$ divides over $B_n$ for every $n < t+1$.
\end{lem}

\noindent
Towards a contradiction, 
\begin{align}\label{assumption leading to a contradiction}
&\text{assume that
$a \in M$ and there are finite
$\es = B_0 \subset B_1 \subset \ldots \subset B_{t+1} \subset M$} \\ 
&\text{such that 
$\tp_\mcM(a / B_{n+1})$ divides over $B_n$ for every $n < t+1$.} \nonumber
\end{align}

\noindent
We will derive a contradiction via a construction of homogeneous simple substructures 
$\mcM = \mcM_0 \supset \mcM_1 \supset \ldots \supset \mcM_{t+1}$ and an argument which is divided
into a few lemmas.

\begin{notation}\label{notation for dividing}{\rm
We will consider dividing in different structures where the universe
of one is included in another. To distinguish which structure we have in mind we use the following notation: if $\mcN$ is a structure,
$\bar{a}, \bar{b} \in N$ and $C \subseteq N$, then $\bar{a} \underset{C}{\ind}^\mcN \bar{b}$ means that 
$\bar{a}$ is independent from $\bar{b}$ over $C$ in $\mcN$ (or ``with respect to $\mcN$'').
}\end{notation}

\noindent
Recall that the (finite) vocabulary of $\mcM$ is denoted $V$.

\begin{defin}\label{definition of R-i}{\rm
Let $p_1(x,y), \ldots, p_t(x,y)$ be an enumeration of $S_2(T)$.
By assumption~(\ref{every 2-type is isolated by an atomic formula}),
there are binary relation symbols $R_1, \ldots, R_t \in V$ such that for each $i$, $R_i(x,y)$ isolates $p_i(x,y)$.
}\end{defin}

\noindent
Note that, since $\mcM$ is homogeneous, every $p \in S_2(T)$ is realized in $\mcM$.

\begin{defin}\label{definition of complete dividing relation}{\rm
Let $\mcN$ be a simple $V$-structure and let $R \in V$ be binary.\\
We call $R$ a {\em dividing relation with respect to $\mcN$} if for all $a, b \in N$, 
$\mcN \models R(a,b)$ implies $a \nind^\mcN b$.
We call $R$ a {\em nondividing relation with respect to $\mcN$} if for all $a, b \in N$, 
$\mcN \models R(a,b)$ implies $a \ind^\mcN b$.
}\end{defin}

\noindent
Note that, in general, a binary $R \in V$ may be neither a dividing relation with respect to $\mcN$ nor a nondividing relation with respect to $\mcN$.

\begin{defin}\label{definition of M-n}{\rm
Let $\mcM_0 = \mcM$ and $T_0 = T$.
For all $n = 1, \ldots, t+1$, let
\begin{align*}
M_n \ &= \ \big\{ a' \in M : \tp_\mcM(a' / B_n) = \tp_\mcM(a / B_n) \big\}, \\
\mcM_n \ &= \ \mcM \uhrc M_n, \ \text{ and}\\
T_n \ &= \ Th(\mcM_n).
\end{align*}
}\end{defin}

\noindent
Hence each $\mcM_n$ is a substructure of $\mcM$ and thus a $V$-structure.
Also note that 
\[ \mcM = \mcM_0 \supset \mcM_1 \supset \ldots \supset \mcM_{t+1} \]
and that $M_t$ is infinite, because $\tp_\mcM(a / B_{t+1})$ divides over $B_t$ and therefore
$\tp_\mcM(a / B_t)$ cannot be algebraic.

\begin{lem}\label{each M-n is simple and homogeneous}
For all $n = 0, \ldots, t$, $\mcM_n$ is simple and homogeneous.
\end{lem}

\noindent
{\bf Proof.}
The case $n = 0$ is trivial, so suppose that $1 \leq n \leq t$.
If a simple structure is expanded with constant symbols (but nothing more) then the resulting
expansion is also simple, by \cite[Remark 2.26]{Cas} for example.
Every infinite structure which is interpretable in a simple structure is simple, by \cite[Corollary 2.8.11]{Wag} for example.
By the $\omega$-categoricity of $\mcM$, $\tp_\mcM(a / B_n)$ is isolated (recall that $B_n$ is finite) and therefore
$M_n$ is $B_n$-definable in $\mcM$. It follows that $\mcM_n$ is interpretable in the expansion of $\mcM$ with constants for elements in $B_n$. 
Thus $\mcM_n$ is simple.

For homogeneity, suppose that $\bar{a} = (a_1, \ldots, a_k), \bar{b} = (b_1, \ldots, b_k) \in (M_n)^k$ and 
$\tp^{at}_{\mcM_n}(\bar{a}) = \tp^{at}_{\mcM_n}(\bar{b})$.
As $\mcM_n$ is a substructure of $\mcM$ we get
\[\tp^{at}_\mcM(\bar{a}) \ = \ \tp^{at}_\mcM(\bar{b}).\]
Since $\bar{a}, \bar{b} \in (M_n)^k$ we also have
\[\tp_\mcM(a_i / B_n) \ = \ \tp_\mcM(b_i /B_n) \ = \ \tp_\mcM(a / B_n) \ \text{ for all $i = 1, \ldots, k$}.\]
Since $\mcM$ is binary we get
$\tp^{at}_\mcM(\bar{a} / B_n) = \tp^{at}_\mcM(\bar{b} / B_n)$ and as $\mcM$ is homogeneous
there is an automorphism $f$ of $\mcM$ such that
$f(\bar{a}) = \bar{b}$ and $f$ fixes $B_n$ pointwise.
Since $M_n$ is $B_n$-definable in $\mcM$, 
it follows that $f$ fixes $M_n$ setwise.
Hence $f \uhrc M_n$ is an automorphism of $\mcM_n = \mcM \uhrc M_n$.
\hfill $\square$

\begin{cor}\label{every R-i isolates a type in every M-n}
For all $i = 1, \ldots, t$ and all $n = 0, 1, \ldots, t$, $R_i(x,y)$ isolates a type in $S_2(T_n)$.
Moreover, every type in $S_2(T_n)$ is isolated by some $R_i$.
\end{cor}

\noindent
{\bf Proof.}
Since $\mcM$ is homogeneous and $\mcM_n \subseteq \mcM$ it follows 
from Definition~\ref{definition of R-i} that for every $i = 1, \ldots, t$ and every atomic $V$-formula
$\varphi(x,y)$, 
\[ \mcM_n \models \forall x,y \big( R_i(x,y) \ \rightarrow \ \varphi(x,y) \big) \ \vee \ 
\forall x,y \big( R_i(x,y) \ \rightarrow \ \neg\varphi(x,y) \big).\]
By Lemma~\ref{each M-n is simple and homogeneous},
$\mcM_n$ is homogeneous and hence it has elimination of quantifiers, 
so the first claim of the corollary follows.
The second claim is just a restatement of what is said in
Definition~\ref{definition of R-i}.
\hfill $\square$
\\

\noindent
Corollary~\ref{every R-i isolates a type in every M-n} immediately implies the following:

\begin{cor}\label{each R-i is a dividing relation or a nondividing relation}
For all $i = 1, \ldots, t$ and all $n = 0, 1, \ldots, t$, $R_i$ is a dividing relation or a nondividing relation (but not both)
with respect to $\mcM_n$.
\end{cor}

\begin{lem}\label{independence in M-n coincides with independence in M over B-n}
For all $n = 0, \ldots, t$ and all $c, d \in M_n$, $c \ind^{\mcM_n} d$ if and only if $c \underset{B_n}{\ind}^{\mcM} d$.
\end{lem}

\noindent
{\bf Proof.}
Since $B_0 = \es$ and $\mcM_0 = \mcM$ the lemma is trivial for $n = 0$.
Let $1 \leq n \leq t$ and $c, d \in M_n$.
Suppose that $c \nind^{\mcM_n} d$.
By the definition of dividing there are a $V$-formula $\varphi(x,y)$ (without parameters)
and $d_i \in M_n$ for $i < \omega$ such that
$\varphi(x,y) \in \tp_{\mcM_n}(c, d)$, 
$\tp_{\mcM_n}(d_i) = \tp_{\mcM_n}(d)$ for all $i$ and
$\{\varphi(x,d_i) : i < \omega\}$ is $k$-inconsistent (with respect to $T_n = Th(\mcM_n)$) for some $k < \omega$.
(It follows from the homogeneity of $\mcM_n$ and
Fact~\ref{facts about simple omega-categorical structures} 
that such $d_i$ can be found in $\mcM_n$.)
Since $\mcM_n$ is homogeneous we may assume that $\varphi$ is quantifier free.
By the definition of $\mcM_n$ we have
$\tp_\mcM(d_i / B_n) \ = \  \tp_\mcM(d / B_n)$ for all $i$.

Since $\mcM_n \subseteq \mcM$, $\varphi(x,y) \in \tp_{\mcM_n}(c,d)$ 
and $\varphi(x,y)$ is quantifier free it follows that
\[ \varphi(x,y) \in \tp_\mcM(c, d). \]
Let $\bar{b}_n$ be an enumeration of $B_n$ and let $\psi(x, y, \bar{b}_n)$ isolate
$\tp_\mcM(c, d / B_n)$. By the homogeneity of $\mcM$ we may assume that $\psi$ is quantifier free.
Then let $\varphi'(x, y, \bar{b}_n)$ be the (quantifier free) formula
\[ \varphi(x, y) \ \wedge \ \psi(x, y, \bar{b}_n), \]
so $\varphi'(x, d, \bar{b}_n) \in \tp_\mcM(c / \{d\} \cup B_n)$.
If 
$\{ \varphi'(x, d_i, \bar{b}_n) : i < \omega\}$
would not be $k$-inconsistent with respect to $\mcM$, then there would be $i_1, \ldots, i_k$ and $c' \in M$ such that
\[\mcM \models \bigwedge_{j=1}^k \varphi(c', d_{i_j}, \bar{b}_n),\]
which, by the definitions of $\mcM_n$ and $\varphi'$,
implies that $c' \in M_n$ and $\mcM_n \models  \bigwedge_{j=1}^k \varphi(c', d_{i_j})$.
Then  $\{\varphi(x, d_i) : i < \omega\}$ is not $k$-inconsistent
with respect to $\mcM_n$, contradicting our assumptions.
Hence $\{ \varphi'(x, d_i, \bar{b}_n) : i < \omega\}$ is $k$-inconsistent 
with respect to $\mcM$ and therefore $c \underset{B_n}{\nind}^\mcM d$.

Now suppose that $c \underset{B_n}{\nind}^\mcM d$ and let $\bar{b}_n$ enumerate $B_n$.
Then there are $\varphi(x, y, \bar{b}_n) \in \tp_\mcM(c, d / B_n)$ and $d_i \in M$ for $i < \omega$
such that $\tp_\mcM(d_i / B_n) = \tp_\mcM(d / B_n)$ for all $i$ and
$\{ \varphi(x, d_i, \bar{b}_n) : i < \omega \}$ is $k$-inconsistent
with respect to $\mcM$ for some $k < \omega$.
By homogeneity of $\mcM$ we may assume that $\varphi$ is quantifier free.
Recall that we assume that $c, d \in M_n$. Since $\tp_\mcM(d_i / B_n) = \tp_\mcM(d / B_n)$ 
we have $d_i \in M_n$ for all $i$.
Since $\mcM_n \subseteq \mcM$ where both structures are homogeneous we also get
\[ \tp_{\mcM_n}(d_i) \ = \ \tp_{\mcM_n}(d) \ \text{ for all $i$}. \]
Let $\psi(x, y)$ isolate $\tp_{\mcM_n}(c, d)$. As $\mcM_n$ is homogeneous we may assume that
$\psi$ is quantifier free.

Suppose for a contradiction that $\{ \psi(x, d_i) : i < \omega \}$ is not $k$-inconsistent with respect to $\mcM_n$.
Then there are $i_1, \ldots, i_k$  and $c' \in M_n$ such that
$\mcM_n \models \bigwedge_{j=1}^k \psi(c', d_{i_j})$.
Because $\psi$ is quantifier free and $\mcM_n \subseteq \mcM$ we get
$\mcM \models \bigwedge_{j=1}^k \psi(c', d_{i_j})$.
Since $\psi$ isolates $\tp_{\mcM_n}(c, d)$ and $c', c, d, d_{i_j} \in M_n$ where $\mcM_n \subseteq \mcM$ it follows
that $\tp^{at}_\mcM(c', d_{i_j}) = \tp^{at}_\mcM(c, d)$ (for all $j$).
By the definition of $M_n$ we have 
$\tp_\mcM(c' / B_n) = \tp_\mcM(c / B_n) = \tp_\mcM(d / B_n) = \tp_\mcM(d_{i_j} / B_n)$.
Since $\mcM$ is binary it follows that 
$\tp^{at}_\mcM(c', d_{i_j} / B_n) = \tp^{at}_\mcM(c, d / B_n)$ which by the homogeneity of $\mcM$
gives 
$\tp_\mcM(c', d_{i_j} / B_n) = \tp_\mcM(c, d / B_n)$ for all $j$.
Consequently $\mcM \models \bigwedge_{j=1}^k \varphi(c', d_{i_j}, \bar{b}_n)$
so $\{ \varphi(x, d_i, \bar{b}_n) : i < \omega \}$ is not $k$-inconsistent
with respect to $\mcM$, contradicting the assumption.
Hence $\{ \psi(x, d_i) : i < \omega \}$ is $k$-inconsistent with respect to $\mcM_n$ 
and therefore $c \nind^{\mcM_n} d$.
\hfill $\square$

\begin{lem}\label{nondividing relation becomes dividing relation}
Let $1 \leq n \leq t$.
There is $1 \leq s \leq t$ such that $R_s$ is a nondividing relation with respect to $\mcM_n$ and, for every $m < n$,
$R_s$ is a dividing relation with respect to $\mcM_m$.
\end{lem}

\noindent
{\bf Proof.}
Recall that for every $n = 0, \ldots, t$, $\mcM_n$ is infinite and (by 
Lemma~\ref{each M-n is simple and homogeneous})  simple and homogeneous.
Fix any $1 \leq n \leq t$.
Let  $\bar{b}_n$ enumerate $B_n$.
For every $0 \leq m < n$ let  $\alpha_m = pw(\bar{b}_n / B_m)$ where the preweight is taken with respect to $\mcM = \mcM_0$.
By Lemma~\ref{the preweight is finite}, $\alpha_m < \omega$ for every $m < n$.
Then let $\alpha = \max\{\alpha_0, \ldots, \alpha_{n-1}\}$, so $\alpha < \omega$.

By repeatedly applying the existence of nondividing extensions and 
Fact~\ref{facts about simple omega-categorical structures} 
to $\mcM_n$ it follows that there are
$c_i \in M_n$ for $i < \omega$ such that $\{ c_i : i < \omega \}$ is an independent set over $\es$ with respect to  $\mcM_n$.
By Ramsey's theorem \cite[Theorem~11.1.3 or its corollary]{Hod}
there are distinct $k_0 < \ldots < k_\alpha < \omega$ such that 
\[ \tp_{\mcM_n}(c_{k_i}, c_{k_j}) \ = \ \tp_{\mcM_n}(c_{k_{i'}}, c_{k_{j'}}) \]
whenever $0 \leq i < j \leq \alpha$ and $0 \leq i' < j' \leq \alpha$.
In other words, by renaming elements for notational simplicity, we have found distinct $d_0, \ldots, d_{\alpha} \in M_n$ such that
$\{d_0, \ldots, d_{\alpha}\}$ is independent over $\es$ with respect to $\mcM_n$ and
$\tp_{\mcM_n}(d_i, d_j) = \tp(d_{i'}, d_{j'})$ whenever $0 \leq i < j \leq \alpha$ and $0 \leq i' < j' \leq \alpha$.
By Corollary~\ref{every R-i isolates a type in every M-n}, 
for some $1 \leq s \leq t$, $R_s$ isolates $\tp_{\mcM_n}(d_i, d_j)$ for all $0 \leq i < j \leq \alpha$.
Hence $\mcM_n \models R_s(d_i, d_j)$  and  thus $\mcM_m \models R_s(d_i, d_j)$ for all $m < n$ and
all $0 \leq i < j \leq \alpha$ (because $\mcM_n \subseteq \mcM_m$ if $m < n$).
As $\mcM = \mcM_0$ we have $\mcM \models  R_s(d_i, d_j)$ for all $0 \leq i < j \leq \alpha$.
Since $\{d_0, \ldots, d_\alpha\}$ is independent over $\es$ in $\mcM_n$ it follows that 
$R_s$ is a nondividing relation with respect to $\mcM_n$.

Let $m < n$ and suppose for a contradiction that $R_s$ is a nondividing relation with respect to $\mcM_m$.
Then $d_i \ind^{\mcM_m} d_j$ for all $0 \leq i < j \leq \alpha$, which by
Lemma~\ref{independence in M-n coincides with independence in M over B-n}
gives $d_i \underset{B_m}{\ind}^{\mcM} d_j$ for all $0 \leq i < j \leq \alpha$.
Now Corollary~\ref{pairwise independence implies independence}
implies that $\{ d_0, \ldots, d_\alpha \}$ is independent over $B_m$ in $\mcM$.
By assumption~(\ref{assumption leading to a contradiction}), the definition of $\mcM_n$ and
since $d_i \in M_n$ and $n > m$ (so $B_m \subseteq B_{n-1}$) we also have $d_i \underset{B_m}{\nind}^\mcM \bar{b}_n$,
for every $i \leq \alpha$.
But as $|\{ d_0, d_1, \ldots, d_\alpha \}| > \alpha \geq \alpha_m$, this contradicts the assumption that 
$pw(\bar{b}_n / B_m) = \alpha_m$.
Hence $R_s$ is (by Corollary~\ref{each R-i is a dividing relation or a nondividing relation}) 
a dividing relation with respect to $\mcM_m$.
\hfill $\square$
\\

\noindent
With Lemma~\ref{nondividing relation becomes dividing relation} we can now derive a
contradiction which proves Lemma~\ref{no long dividing sequence of real elements in M} 
and hence also Theorem~\ref{binary homogeneous simple structures have finite rank}.
By Lemma~\ref{nondividing relation becomes dividing relation},
there is $1 \leq s \leq t$ such that $R_s$ is a nondividing relation with respect to $\mcM_t$
and a dividing relation with respect to $\mcM_n$ for every $n < t$. 
Without loss of generality (by just reordering $R_1, \ldots, R_t$ if necessary) we can assume that $s = 1$.
By Lemma~\ref{nondividing relation becomes dividing relation} again,
there is $1 \leq s \leq t$ such that $R_s$ is a nondividing relation with respect to $\mcM_{t-1}$
and a dividing relation with respect to $\mcM_n$ for every $n < t-1$.
By the previous step we must have $s > 1$. Without loss of generality
(by just reordering $R_2, \ldots, R_t$ if necessary) we may assume that $s = 2$. 
If we continue in the same way until $t$ steps are finished we find that each one of
$R_1, \ldots, R_t$ is a dividing relation with respect to $\mcM_0 = \mcM$. 
Since $\mcM$ is simple there are (by the existence of nondividing extensions)
$c, d \in M$ such that $c \ind^\mcM d$. Then $\tp_\mcM(c,d)$ is isolated by
$R_i$ for some $1 \leq i \leq t$ and this $R_i$ must be a nondividing relation with respect to $\mcM$,
contradicting (via Corollary~\ref{each R-i is a dividing relation or a nondividing relation}) 
that every $R_i$ is a diving relation with respect to $\mcM$.
This finishes the proof of Lemma~\ref{no long dividing sequence of real elements in M} 
and of Theorem~\ref{binary homogeneous simple structures have finite rank}.

The assumption that $\mcM$ is binary was used directly in the proofs of
Lemmas~\ref{each M-n is simple and homogeneous}
and~\ref{independence in M-n coincides with independence in M over B-n}
and indirectly in the proof of Lemma~\ref{nondividing relation becomes dividing relation}
through the use of 
Corollary~\ref{pairwise independence implies independence}.

\end{document}